\documentclass [a4,10pt,twoside]{article}
\usepackage{amsmath, amssymb, amsthm, amsfonts, amsxtra, latexsym, amscd,
pb-diagram,  graphics,rotating,xy}

\theoremstyle{plain}

\newcommand{\A}{\mathcal{A}}

\newcommand{\R}{\mathcal{R}}
\newcommand{\tx}{\otimes}
\newcommand{\ts}{\oplus}
\begin{document}
\title{ON THE CATEGORICAL INTERPRETATION OF THE RING COHOMOLOGY}
\author{N.T.Quang and N. T. Thuy}
\pagestyle{myheadings} 
\markboth{On the categorical interpretation of the ring cohomology}{N. T. Quang and N. T. Thuy}
\maketitle
\setcounter{tocdepth}{1}
%\tableofcontents 
As known, for the ring structure, the two properties: association and commutation respect to +, can be replaced with the  ``association-commutation'' property:
$$(a+b)+(c+d)=(a+c)+(b+d).$$ 
This leads to two ways to construct the ring cohomology: cohomology for rings which are regarded as algebras due to Shukla [Sh], and ring cohomology due to Mac Lane [Mac].

The ring structure has been categorized by the definition of an Ann-category [4]. Each congruence class of Ann-categories is characterized by a cohomology class of structures (Theorem 4.3[4]). For regular Ann-categories, this class is bijective to the group $H^3_{Sh}(R,M)$ (Theorem 4.4 [5]). In the general case, this group is replaced with the group $H^3_{Mac}(R,M)$ (Theorem 7.6 [7]).

In [1], authors have modified the definition of an Ann-category to be the one of a {\it categorical ring,} where the condition {\it(Ann-1)} is omitted, and the compatibility of the operation $\tx$ with the associativity and commutativity constraints respect to + is replaced with the compatibility of the operation $\tx$ with the ``associativity - commutativity'' constraint. Categorical rings have been classified by the Mac Lane cohomology group $H^{3}_{Mac}(R, M).$

There are 2 problems emerged: Firstly, do the sets of categorical rings and Ann-categories coincide? Secondly, in the case that these sets can not be proved to coincide, is it possible to prove the existence of the bijection between the set of congruence classes of pre-sticked Ann-categories of the type $(R,M)$ and the group $H^{3}_{Mac}(R, M)?$ We have completely solved both of these problems.

In [6], authors have showed that the set of Ann-categories is a subset of the set of categorical rings, and these two sets coincide if and only if the operation $\tx$ in each categorical ring is compatible with the unitivity constraints of the operation $\ts.$ (This condition is also called the condition $(U)$). A question is raised: May the condition $(U)$ be deduced from the axiomatics of a categorical ring? The second problem has been solved in [7].

After that, in [2], authors have proved that the condition $(U)$ is always satisfied in a categorical ring, and therefore the two sets of the categories coincide. However, authors have made a mistake in the proof of the Lemma.

Indeed, in [6], in order to prove that a categorical ring together with the condition $(U)$ is an Ann-category, it is essential to prove that the operation $\tx$ is compatible with the {\it associativity} constraint respect to the operation $\ts.$ Therefore, the Proposition referred in [2] is unable to be applied.

We now show more clearly some mistakes in [1].

1. If $\A$ is an Ann-category, then the set $\pi_0(\A)$ of congruence classes of objects is a ring with the operation induced by the operation $\tx,$ and $\Pi_1(\A)=Aut(0)$ is a $\Pi_0(\A)$-bimodule with the actions given by:
$$su=\lambda_X(u),\ \ us=\rho_X(u),\ \ X\in s\in\Pi_0(\A),u\in\Pi_1(\A)$$
This fact is possible thanks to the existence of isomorphisms (see Proposition 1.4 [5]):
$$A\otimes 0\rightarrow 0 \ ;\ 0\otimes A\rightarrow 0.$$
Unfortunately, for categorical rings, the fact that these isomorphisms exist is still an open question, so we even can not assert that $\Pi_1(\A)=Aut(0)$ has the structure of a $\Pi_0(\A)$-bimodule!

2. The mistake of using the property $(U)$ repeatedly occurs in the proof of Proposition 3.1, at the beginning of Section 4. Since the isomorphisms $\overline{r}\otimes 0\rightarrow 0 \ ;\ 0\otimes \overline{s}\rightarrow 0$ do not exist, the determination of the morphism $\varphi_{\bullet}(r,s,t)\in\pi_1$ from the morphism $\overline{rst}\rightarrow \overline{rst}$ in $\mathcal R$ is illegal. It is similar for $\varphi_{\bullet +}, \varphi_{+ \bullet }, \varphi_{+}.$

3. In Section 4 [1], afer defining the collection of functions $$\varphi=(\varphi_{\bullet}, \varphi_{\bullet +}, \varphi_{+ \bullet }, \varphi_{+}),$$ authors have proved that $\varphi$ is a 3-cocycle of Mac Lane cohomology. We can see that the isomorphisms $\sigma_{\bullet}(r,s), \sigma_{+}(r,s)$ referred in this section turn out to be a {\it stick,} and the way to determine $\varphi_{\bullet +}, \varphi_{+ \bullet }$ turns out to be the way to determine {\it the constraints} $\alpha, \lambda, \rho$ of a reduced Ann-category in [4,5,7].

The most important proof in this section is to use 8 ``reduced'' diagrams to prove that the collection $\varphi$ is a 3-cocycle. In [7], this has been proved in another way, avoiding the complexity of diagrams. So the functions $\alpha, \lambda, \rho$ together with $\xi,\eta$ (respect to the associativity, commutativity constraints of + in the reduced Ann-category) is a structure of the reduced Ann-category $\mathcal S.$ Therefore, the fact that the collection $f=(\xi,\eta,\alpha, \lambda, \rho)$ satisfies the relations referred in Theorem 3.1 [5] follows from the definition of an Ann-category. We now replace the 2 functions $\xi,\eta$ with the function $\sigma$ by:
$$\sigma(x,y,z,t)= \xi(x+y,z,t)-\xi(x,y,z)+\eta(y,z)+\xi(x,z,y)-\xi(x+z,y,t) $$
Then, the relations in Theorem 3.1 [5] will give us the relations of a 3-cocycle (Proposition 7.3 [7]). Also, we note that, while in [1] there is no information of the relation between the structure of $\R$ and the structure of $ \mathcal S=(\Pi_{0}, \Pi_{1}),$ in [4,5,7] we assert that they are Ann-equivalent.

4. Some results in [1] rely on the remark: For any two parallel morphisms $\varphi, \varphi':x\rightarrow y,$ there exists a {\it unique} $\alpha\in \Pi_{1}(\A)$ making the diagram
\[
\begin{diagram}
\node{0+x} \arrow{e,t}{\alpha+\varphi}\arrow{s,l}{\lambda(x)}\node{0+y} \arrow{s,r}{\lambda(y)}\\ 
\node{x} \arrow{e,t}{\varphi'}\node{y}
\end{diagram}
\]
commute. This remark has not been proved so its correctness can not be verified.

5. The morphism $$\left\langle \begin{array}{cccc}
a&b\\
c&d
\end{array}\right\rangle:(a+b)+(c+d)\to(a+c)+(b+d)$$ 
in [1] is just the morphism $v$ in [4,5,7]. Since $(\mathcal A,\oplus)$ is a symmetric monoidal category, the morphism $v$ built from $a_+, c_+$ is unique thanks to the coherence theorem. That means every way from $(a+b)+(c+d)$ to $(a+c)+(b+d)$ is the same. So the diagram in [1] can be replaced with the diagram
\[
\begin{diagram}
\node{(x+y)+(z+t)}\arrow{s,l}{v} \arrow{e,t}{a_{+}}\node{((x+ y)+ z)+ t}\node{(x+ (y+ z))+ t}\arrow{w,t}{a_{+}+ t}\arrow{s,r}{(x+z)+ t}\\
\node{(x+z)+(y+t)}\arrow{e,t}{a_{+}}\node{((x+z)+y)+t}\node{(x+(z+y))+ t.}\arrow{w,t}{a_{+}+t}
%\tag{5}
\end{diagram}
\]

6. In proposition 4.3 [1], the assumption: ``there exists a 2-homomorphism ${\bf f}:\mathcal R\rightarrow \mathcal R'$ inducing identity maps on $\Pi_0$ and $\Pi_1$'' is so strict. We have known that if ${\bf f}$ is an equivalence, then it induces the ring isomorphism $F_0:\Pi_0\rightarrow \Pi_0'$ and the group isomorphism $F_1:\Pi_1\rightarrow \Pi_1'$ satisfying the relations:
$$F_0(ru)=F_0(r)F_1(u)\ , \ F_0(us)=F_1(u)F_0(s)$$
Therefore, the definition of a {\it pre-stick} as Definition 4.1 [5] is necessary to fix the pair $(B,R).$

7. We now consider the construction of the inverse map
$$H^{3}(B,R)\rightarrow Crext(R,B).$$
After defining the constraints on the category $R_{\varphi},$ the proof of the fact that $R_{\varphi}$ is a categorical ring is not completed, because there is no relation showing the compatibility of the associativity constraint with the commutativity one respect to +. This has been showed in the proof of Proposition 7.3 [7].

8. In order to prove Theorem 4.4, authors of [1] built a 2-homomorphism $f:R_{\varphi}\rightarrow R.$ This turns out to be the Ann-functor $H:\mathcal S\rightarrow \mathcal R,$ and the choice of the functions:
$$\sigma_+(0,r)=\lambda(\overline{r})\ ,\ \sigma_+(r,0)=\rho(\overline{r}) $$
is just the choice of the {\it sticks} $\psi_{0,r}\,\ \psi_{r,0} $ in [4,5]. However, without the condition $(U)$, the functions $\lambda(\overline{r}), \rho(\overline{r})$ have not been determined, i.e., ${\bf f}$ can not be determined.

Finally, the fact that ${\bf f}$ is a 2-homomorphism is proved uneasily (see Section 2 [8]). In [5], ${\bf f}$ turns out to be the Ann-functor $H.$ In our opinion, most of the proofs in [1], fundamentally, are similar to the ones in [5], just different in logical order and the detail of interpretaion.
\begin{center}

\end{center}
\end{document}